\documentclass[a4paper,12pt,reqno]{amsart}
\usepackage{amsmath}
\usepackage{amssymb}
\usepackage{amsthm}

\usepackage{amscd}
\usepackage{amsfonts}
\usepackage{setspace}
\usepackage{version}

\title[A more intuitive proof of a sharp version of Hal\'{a}sz's theorem]{A more intuitive proof of a sharp version of Hal\'{a}sz's theorem}
\author{Andrew Granville}
\address{AG: D\'epartement de math\'ematiques et de statistique\\
Universit\'e de Montr\'eal\\
CP 6128 succ.~Centre-Ville\\
Montr\'eal, QC H3C 3J7\\
Canada; 
and Department of Mathematics \\
University College London \\
Gower Street \\
London WC1E 6BT \\
England.
}
  \email{{\tt andrew@dms.umontreal.ca}}

  \author{Adam J Harper}
\address{AJH: Mathematics Institute, Zeeman Building, University of Warwick, Coventry CV4 7AL, England}
\email{\tt A.Harper@warwick.ac.uk}
\author{K. Soundararajan}
\address{KS: Department of Mathematics, Stanford University, Stanford CA 94305, USA}
\email{{\tt  ksound@stanford.edu}}
\date{12th June 2017}
\thanks{Andrew Granville has received funding in aid of this research from the
European Research Council  grant agreement n$^{\text{o}}$ 670239, and from NSERC Canada under the CRC program. Adam Harper was supported, for parts of the research, by a postdoctoral fellowship from the Centre de recherches math\'{e}matiques in Montr\'{e}al, and by a research fellowship at Jesus College, Cambridge. Kannan Soundararajan was partially supported by NSF grant DMS 1500237, and a Simons Investigator grant from the Simons Foundation.}

\numberwithin{equation}{section}
\theoremstyle{plain}

\onehalfspace
\newcommand{\N}{\mathbb{N}}

\newcommand{\C}{\mathbb{C}}

\newtheorem{thm1}{Theorem}
\newtheorem{lem1}{Lemma}

\theoremstyle{remark}
\newtheorem{remark}[subsection]{Remark}

\begin{document}

\maketitle

\begin{abstract}
We prove a sharp version of Hal\'{a}sz's theorem on sums $\sum_{n \leq x} f(n)$ of multiplicative functions $f$ with $|f(n)|\le 1$. Our proof avoids the ``average of averages'' and ``integration over $\alpha$'' manoeuvres that are present in many of the existing arguments. Instead, motivated by the circle method we express $\sum_{n \leq x} f(n)$ as a triple Dirichlet convolution, and apply Perron's formula.
\end{abstract}

\section{Introduction}
Given a multiplicative function $f : \N \rightarrow \C$, for each $x \geq 2$ let its summatory function be
\[
S(x) := \sum_{n \leq x} f(n), \text{ and  }   F_{x}(s) := \prod_{p \leq x} \left(1 + \sum_{k=1}^{\infty} \frac{f(p^{k})}{p^{ks}} \right)
\]
 denote the corresponding truncated Euler product.

In this note we shall prove the following form of Hal\'{a}sz's theorem on mean values of multiplicative functions taking values in the 
unit disc.

\begin{thm1}
For $f, F_x$ and $S(\cdot)$   as above, suppose that $|f(n)| \leq 1$ for all integers $n\geq 1$.  Define the quantity $L(x)$ by setting 
$$ 
L(x)^2 :=  \sum_{|N| \leq \log^{2}x + 1} \frac{1}{N^{2}+1} \sup_{|t-N| \leq 1/2} |F_{x}(1+it)|^{2}  . 
$$
Then we have
$$
 |S(x)| \ll x \frac{L(x)}{\log x} \log\Big(100\frac{\log x }{L(x)}\Big) + x \frac{\log\log x}{\log x}  . 
 $$
\end{thm1}

This is essentially the same as the version of Hal\'{a}sz's theorem proved by Montgomery~\cite{montmult} (if we note that $|F_{x}(1+it)|^{2} \asymp |F(1+\frac 1{\log x} + it)|^2$, where $F(s)$ denotes the full Euler product over all primes), and is known to be quantitatively sharp (see the papers of Granville and Soundararajan~\cite{gransoundmean} and Montgomery~\cite{montmult}). See Hal\'asz's papers~\cite{Hal1, Hal2} for his original arguments which were refined by Montgomery~\cite{montmult}, and see Chapter III.4 of Tenenbaum~\cite{tenenbaumintro} for a textbook treatment.

\vspace{12pt}
Our proof here is, hopefully, more intuitive and easier to motivate than the existing proofs, although it also has important features in common with several of them. We begin by expressing $S(x)$ as a triple Dirichlet convolution, and using Perron's formula to relate our triple Dirichlet convolution  to the Dirichlet series $F_{x}(s)$ and two other Dirichlet polynomials. This is done by analogy with the circle method, as we want to use a pointwise bound for $F_{x}(s)$ and obtain a mean square bound for the remaining two Dirichlet polynomials. To carry everything out with little loss, we break the Dirichlet convolution into subsums which depend on the size of one of the variables $p$.   Our proof avoids the ``average of averages'' step in many other treatments of Hal\'{a}sz's theorem, and in particular it avoids the arguably slightly obscure ``integration over $\alpha$'' device from many of the treatments.

Our   longer companion paper~\cite{granharpsound}  uses a similar strategy  to prove various generalisations of Hal\'{a}sz's theorem,  including for multiplicative functions bounded by divisor functions, and treating sums over short intervals and arithmetic progressions.  However we give here the original argument, stripped of   the   technicalities in the more general argument of
~\cite{granharpsound} (compare, for example, the more complicated but more easily generalisable triple convolution in ~\cite{granharpsound}).

\section{A lemma concerning prime numbers}
We will need some basic information about the integrals of Dirichlet polynomials supported on the primes. We record a suitable result here.

\begin{lem1}
Uniformly for any complex numbers $(a_n)_{n=1}^{\infty}$ and any $T \geq 1$, we have
$$ \int_{-T}^{T} \Biggl|\sum_{T^{2} \leq n \leq x} \frac{a_n \Lambda(n)}{n^{1+it}} \Biggr|^2 dt \ll \sum_{T^{2} \leq n \leq x} \frac{|a_n|^2 \Lambda(n)}{n} . $$
\end{lem1}

\begin{proof}[Proof of Lemma 1]
This follows by inserting a smooth weight $\Phi(t/T)$ into the integral, expanding out, and applying a Brun--Titchmarsh upper bound for primes in short intervals at a suitable point. See Lemma 2.6 of \cite{granharpsound}, for example, for a full proof. 
\end{proof}

Mean value results and majorant principles of this kind  are often used in multiplicative number theory, and proved in very similar ways (see, e.g., the Lemma in section 2 of Montgomery~\cite{montmult}); some would be sufficient for our purposes. However the most standard mean value theorem for Dirichlet polynomials, which implies that $\int_{-T}^{T} |\sum_{n \leq x} \frac{a_n}{n^{it}} |^2 dt = \sum_{n \leq x} |a_n|^2 (2T + O(n))$, would {\em not} suffice because in Lemma 1 it would yield a multiplier $\Lambda(n)^2$, rather than $\Lambda(n)$, on the right hand side.

\section{Proof of Theorem 1: the combinatorial part}
We begin  by expressing $S(x)$ as a triple Dirichlet convolution, up to an acceptable error. Since $\log n= \sum_{d|n} \Lambda(d)$ we have 
$$ 
\sum_{n\le x} f(n) \log n = \sum_{d\le x} \Lambda(d)\sum_{m\le x/d} f(md) = \sum_{mp\le x} f(m) f(p) \log p + O(x), 
$$
the error term arising from bounding trivially the contribution of prime power values of $d$, and the terms with $(m,d)>1$.   Since 
$$ 
\sum_{n\le x} \log (x/n) = O(x),
$$ 
we deduce that 
$$ 
S(x) = \frac{1}{\log x} \sum_{n\le x} f(n) (\log n + \log x/n) = \frac{1}{\log x}\sum_{mp\le x} f(m) f(p) \log p  + O\Big(\frac{x}{\log x}\Big). 
$$ 
This is a double multiplicative convolution, since we have the two variables $p$ and $m$ in the sum.

We  repeat the above argument  to arrive at a triple convolution.  For technical convenience we begin by discarding those primes $p$ for which $p \leq \log^{4}x$ or $p > x/2$ from the sum, which gives rise to an acceptable error term $O(x \frac{\log\log x}{\log x})$.  For primes $p$ in the range $\log^{4}x<p\leq x/2$, we use the above argument to replace the sum over $m\le x/p$ by a double convolution; that is,
$$ 
\sum_{m\le x/p} f(m) = S(x/p) =\frac{1}{\log (x/p)}  \sum_{nq\le x/p}  f(n) f(q) \log q  + O\Big(\frac{x}{p\log x/p}\Big).  
$$ 
Therefore 
\begin{align} 
\label{3.1}
S(x) &= \frac{1}{\log x} \sum_{\log^{4}x<p\leq \frac x2} f(p)\log p \sum_{m\le x/p} f(m) + O\Big( x\frac{\log \log x}{\log x}\Big) 
\nonumber \\
&= \frac{1}{\log x} \sum_{\log^{4}x<p\leq \frac x2} \frac{ f(p)\log p}{\log (x/p)} \sum_{nq\le x/p}  f(n) f(q) \log q+ O\Big(x \frac{\log \log x}{\log x}\Big), 
\end{align}
since 
$$ 
\sum_{\log^{4}x<p\leq \frac x2} \frac{\log p}{p\log (x/p)} \ll \log \log x. 
$$ 
 We have arrived at the desired triple multiplicative convolution. 
 


\vspace{12pt}

The range of $q$ in (3.1) is severely restricted when $p$ is large, which will lead to bigger  error terms, so it pays to treat summands differently depending on the size of $p$. We achieve this by partitioning up the range of the primes $p \in {\mathcal P}= [(\log x)^4,x/2]$ into the intervals ${\mathcal P}_k = {\mathcal P} \cap (x^{1-e^{1-k}},x^{1-e^{-k}}]$, where $k$ runs through the integers from $1$ to $\log \log x + O(1)$.   Define 
$$ 
S_k(x) = \sum_{\substack{pqn\le x \\ p\in {\mathcal P_k}} } \frac{f(p)\log p}{\log (x/p)} f(n) f(q) \log q , 
$$ 
so that (3.1) implies
$$ 
S(x) \ll \frac{1}{\log x} \sum_{k=1}^{\log \log x+O(1)} |S_k(x)| + x\frac{\log \log x}{\log x}.
$$ 


Since each $|f(p)f(q)f(n)| \le 1$, we may bound $S_k(x)$ trivially as follows:  
\begin{eqnarray}
|S_k(x)| & \leq & \sum_{\substack{pq \le x \\ p\in {\mathcal P_k}} } \frac{ \log p}{\log (x/p)}  \log q  \sum_{n\leq x/pq} 1  \nonumber \\
& \leq & x \sum_{\substack{  p\in {\mathcal P_k}} }   \frac{\log p}{p \log(x/p)} \sum_{q\leq x/p} \frac{\log q}{q}   \nonumber \\
& \ll & x \sum_{x^{1-e^{1-k}} < p \leq x^{1-e^{-k}}} \frac{\log p}{p} \ll e^{-k} x \log x.  \nonumber
\end{eqnarray}
Thus the sum of $|S_k(x)|$ over all integers $k > \log(100\log x /L(x))$ (where $L(x)$ is as in the statement of Theorem 1) 
leads to a bound that is acceptable for Theorem 1.   
 
To complete the proof of  Theorem 1, it therefore suffices to show that 
\begin{equation}\label{remainstoprove}
S_k(x) \ll x L(x) + x 
\end{equation}
for all positive integers $k \leq \log(100\log x /L(x))$.

\vspace{12pt}
\begin{remark}
Partitioning the range for $p$ and applying the triangle inequality might appear to be wasteful. However, in the worst case, there is no loss in introducing absolute values, since the arguments of the values of the $f(p)$ with $p > x^{1-e^{-1}}$ could have been chosen, given the values $f(q^k)$ for $q \leq x^{1-e^{-1}}$, so that   $f(p)$ times the sum over $qn\leq x/p$    all point in exactly the same direction. Indeed,  extremal examples for Hal\'{a}sz's theorem can behave precisely in this way, as in the introduction to Montgomery's paper~\cite{montmult}.
\end{remark}

\begin{remark}
If the multiplicative function $f(n)$ is supported only on numbers with all their prime factors $\leq x^{0.999}$, say, (that is, the $x^{0.999}$-smooth numbers), then there will only be a bounded number of terms $k$ in our decomposition of the $p$-sum.  For such functions $f$, Hal\'{a}sz's theorem can be improved, using (3.2), to
$$
 |S(x)| \ll   \frac{x}{\log x}   (L(x)+1) ,
 $$
after taking a little more care in handling the discarded contribution from the primes $p \leq \log^{4}x$. As far as we know, this has not been noted previously.
\end{remark}

\section{Proof of Theorem 1: the analytic part}
It remains  to prove \eqref{remainstoprove}.  If $p$ lies in ${\mathcal P_k}$ then, since $pq$ must be at most $x$, the 
prime $q$ is constrained to $q \le x^{e^{1-k}}$.   Therefore, using a truncated Perron formula, we get 
$$ 
S_k(x)  \ll x\int_{1-i(\log x)^2}^{1+i (\log x)^2} \Big| \sum_{p \in {\mathcal P}_k} \frac{f(p)\log p}{p^s\log (x/p)} \Big| 
\Big| \sum_{q\le x^{e^{1-k}} } \frac{f(q)}{q^s} \log q \Big| |F_x(s)| \frac{|ds|}{|s|}  + x. 
$$ 
Using the Cauchy--Schwarz inequality, we bound the integral above by $\sqrt{I_1 I_2}$, where 
$$ 
I_1 = \int_{1-i(\log x)^2}^{1+i(\log x)^2} \Big| \sum_{p\in {\mathcal P}_k} \frac{f(p)\log p}{p^s \log (x/p)} \Big|^2 |ds|, 
$$
and 
$$ 
I_2 = \int_{1-i(\log x)^2}^{1+i(\log x)^2} \Big| \sum_{q\le x^{e^{1-k}}} \frac{f(q)}{q^s} \log q \Big|^2 |F_x(s)|^2 \frac{|ds|}{|s|^2}. 
$$ 


Splitting the integral in $I_2$ into intervals of length $1$, we may bound it 
by 
$$ 
I_2 \ll \sum_{|h| \leq \log^{2}x + 1} \frac{1}{h^{2}+1} \sup_{|t-h| \leq 1/2} 
|F_{x}(1+it)|^{2} \int_{1+i(h-1/2)}^{1+i(h+1/2)} \Big| \sum_{q \leq x^{e^{1-k}}} \frac{f(q)}{q^{s}} \log q \Big|^{2} |ds| . 
$$
Recalling that $q$ runs over primes, we can apply Lemma 1 with $T=1$, $a_q = f(q)q^{-ih}$ for primes $q$, and $a_q = 0$ otherwise, and 
deduce that 
\begin{eqnarray}
I_2 & \ll & \sum_{|h| \leq \log^{2}x + 1} \frac{1}{h^{2}+1} \sup_{|t-h| \leq 1/2} |F_{x}(1+it)|^{2} \sum_{q \leq x^{e^{1-k}}} \frac{\log q}{q} \nonumber \\
& \ll &  L(x)^{2} \,    e^{-k} \log x . \nonumber
\end{eqnarray}

To bound $I_1$, we use Lemma 1 again (noting that ${\mathcal P}_k$ only has primes larger than $(\log x)^4$ for all $k$), to 
obtain
$$ 
I_1 \ll \sum_{p\in {\mathcal P}_k}  \frac{\log p}{p \log^{2}(x/p)} \ll \frac{e^{2k}}{\log^{2}x} \sum_{x^{1-e^{1-k}} < p \leq x^{1-e^{-k}}} \frac{\log p}{p} \ll \frac{e^k}{\log x} . $$

Combining the foregoing estimates, we obtain \eqref{remainstoprove} and therefore the bound claimed in Theorem 1.
\qed

\section*{Acknowledgements}
We would like to thank Ben Green for his insistence that ``there must be a more motivated proof'' of Hal\'{a}sz's theorem, which led to this article as well as to \cite{granharpsound}.

\end{document}